\documentclass[12 pt,twoside]{amsart}

\usepackage{amsopn}
\usepackage{amssymb}
\usepackage{amscd}

\newtheorem{theorem}{Theorem}[section]
\newtheorem{lemma}[theorem]{Lemma}
\newtheorem{proposition}[theorem]{Proposition}
\newtheorem{corollary}[theorem]{Corollary}

\theoremstyle{definition}
\newtheorem{definition}[theorem]{Definition}

\theoremstyle{remark}
\newtheorem{remark}[theorem]{Remark}

\numberwithin{equation}{section}

\begin{document}

\title[On the existence of J-class operators on Banach spaces]{On the existence of J-class operators on Banach spaces}

\author[Amir Bahman Nasseri]{Amir Bahman Nasseri}
\address{Fakult\"{a}t f\"{u}r Mathematik, Technische Universit\"{a}t Dortmund, D-44221 Dortmund, Germany}
\email{amirbahman@hotmail.de}
\thanks{During this research the author was fully supported by Prof. Dr. Rainer Brück at the Technische Universität Dortmund, Germany. He would also like to express his gratitude to Prof. Dr. W. Kaballo for some helpful comments.}

\date{}

\begin{abstract}
In this note we answer in the negative the question raised by G.Costakis and A.Manoussos, whether there exists a J-class operator on every non-separable Banach space. In particular we show that there exists a non-separable Banach space constructed by A.Arvanitakis, S.Argyros and A.Tolias such that the J-set of every operator on this space has empty interior for each non-zero vector. On the other hand, on non-separable spaces which are reflexive there always exist a J-class operator.
\end{abstract}

\maketitle

\section{Preliminaries and the main result}
Let $X$ be a real or complex Banach space. If $X$ is a real Banach space then $X_{\mathbb{C}}$ is the complexification of $X$. By $L(X)$ we mean the space of all bounded linear operators acting on $X$. If $T\in L(X)$ the symbol $\sigma(T)$ stands for the spectrum of $T$. Consider any subset $C$ of $X$. The symbol $\stackrel{\circ}{C}$ denotes the interior of $C$ in the norm-topology of $X$. The symbol $\text{orb}(T,x)$ denotes the orbit of $x$ under $T$, i.e. $\text{orb}(T,x):=\{T^nx: n\in\mathbb{N}\}$. If $X$ is separable and $\text{orb}(T,x)$ is dense, then $T$ is called hypercyclic, which is equivalent that $T$ is topologically transitive, i.e. for each non-empty pairs of open subsets $U,V\subset X$, there exists a positve integer $n$, such that $T^{n}(U)\cap V\neq\emptyset$. By $J_{T}(x)$ we mean the J-set of $x$ under $T$, i.e.
\begin{align*}J_T(x):=\{y\in X: \text{ there exists a strictly increasing sequence}&\\ \text{of natural numbers }(k_n) \text{ and a sequence }&\\ (x_n)\text{ in } X, \text{ such that } x_n\rightarrow x \text{ and } T^{k_n}x_n\rightarrow y\}.
\end{align*}
If $J_{T}(x)=X$ for some $x\in X\backslash\{0\}$, then $T$ is called a J-class operator and if $J_{T}(0)=X$ and $J_{T}(x)\neq X$ for every $x\in X\backslash\{0\}$, then $T$ is called a zero J-class operator. By $A_{T}$ we denote the set of all $x\in X$, such that $J_T(x)=X$. On separable spaces every hypercyclic operator is J-class, but the converse is not true. It is shown that on $l^{\infty}$, there does not exist a topological transitive operator. On the other hand there exist J-class operators like the weighted backward shift $\lambda B:l^{\infty}\rightarrow l^{\infty}$, $\lambda B(x_1,x_2,\ldots):=(\lambda x_2,\lambda x_3,\ldots )$ for $\left|\lambda\right|>1$. Therefore it is natural to ask, if there exists always a J-class operator on non-separable Banach spaces.  
Our main result is the following:
\begin{theorem}\label{main}
There exists a non-separable complex Banach space $X$ on which the J-set of every operator has empty interior for every non-zero vector. Consequently there exists no J-class operator on $X$.
\end{theorem}

The space we consider in Theorem 1.1 is a non-separable HI( hereditarily indecomposable) Banach space constructed by Argyros, Arvanitakis and Tolias, \cite{Argyros}. Every bounded operator on this space has the form $T=\lambda I+S$, where $\lambda\in\mathbb{C}$ and $S\in L(X)$ is strictly singular. At this point we begin with some properties of strictly singular operators.

\begin{definition}
Let $X,Y$ be infinite dimensional Banach spaces. A linear and bounded operator $S:X\rightarrow Y$ is called strictly singular, if for every  infinite dimensional subspace $M\subset X$ the restriction $S_{|M}: M\rightarrow S(M)$ is not an isomorphism (linear homeomorphism). 
\end{definition}
\begin{remark} If $X=Y$, then an immediate consequence of the above definition is that $0\in\sigma(S)$.
\end{remark}
The next three theorems can be found in \cite{Abramovich}.
\begin{theorem}
Assume that $S\in L(X)$ is strictly singular and that an
operator $T\in L(X)$ has an at most countable spectrum. Then the spectrum
of S+T is at most countable and zero and the points of $\sigma(T)$ are
the only possible accumulation points of $\sigma(S + T)$.
\end{theorem}
\begin{remark}
In the case that $T=0$ in Thm.1.4, we conclude with Rem. 1.3 that $0\in\sigma(T)=\sigma(S)$, with $0$ as the only possible accumulation point. Therefore 
the spectrum of the operator $\lambda I+S$ is equal to $\lambda+\sigma(S)$ and is at most countable with $\lambda$ as the only possible accumulation point.
\end{remark} 
\begin{theorem}
Let $X$ be a Banach space. The collection of all strictly singular operators is a closed subspace in $L(X)$, which is also a two-sided ideal.
\end{theorem}

\begin{proposition} Let $X$ be a real Banach space and $X_{\mathbb{C}}$ its complexification.
A bounded operator $T:X\rightarrow X$ is strictly singular if and only if $T_{\mathbb{C}}: X_{\mathbb{C}}\rightarrow X_{\mathbb{C}}$ is strictly singular.
\end{proposition}

We show now that a certain class of operators on a complex Banach space can not be J-class. For this purpose we need the following lemma , see \cite{G.C.}.
\begin{lemma}
Let $X$ be a Banach space and $T\in L(X)$. Suppose $\stackrel{\circ}{J_{T}(x)}\neq\emptyset$ for some $x\in X\backslash\{0\}$. Then $\sigma (T)\cap\partial \mathbb{D}\neq\emptyset$.
\end{lemma}
\begin{theorem}
Let $X$ be a complex infinite dimensional Banach space. Then the interior of the J-set $J_{T}(x)$ of the operator $T:X\rightarrow X$ defined by $T:=\lambda id +S$, where $S$ is strictly singular, is empty for all $\left|\lambda\right| >1$ and each non-zero vector $x$. Furthermore $T$ is zero J-class if and only if
$\sigma(T)\subset\mathbb{C}\backslash\overline{\mathbb{D}}$.
\end{theorem}
\begin{proof}
By Rem. 1.5 it follows that $\lambda\in\sigma(T)$ and it is the only possible accumulation point. We decompose the spectrum in $\sigma_1:=\{\mu\in\sigma(T):\left|\mu\right|\leq1\}$ and $\sigma_2:=\{\mu\in\sigma(T):\left|\mu\right|>1\}$. Clearly then $\lambda\in\sigma_2$.
$\sigma_1$ is closed and since $\lambda\in\sigma_2$, and $\lambda$ is the only possible accumulation point, we conclude that $\sigma_2$ is also closed. Futhermore $\sigma_1$ and $\sigma_2$ are disjoint. By the Riez-decomposition theorem we can decompose $X=M_1\oplus M_2$, where $M_1$ and $M_2$ are closed and $T$-invariant subspaces and $\sigma_1=\sigma(T_{|M_1})$,
$\sigma_2=\sigma(T_{|M_2})$. Now assume that there exist a non-zero vector $x$, such that $J_{T}(x)$ has non-empty interior. Then $x=x_1+x_2$, where $x_1\in M_1$ and $x_2\in M_2$. Then since $x$ is not zero, $x_1$ or $x_2$ is not zero. Now we conclude that $M_1$ is finite dimensional. Otherwise suppose $M_1$ is infnite dimensional, then
$S_{|M_1}$ is strictly singular and it follows from Rem. 1.5, that $0\in\sigma(S_{|M_1})$ and therefore $\lambda\in\sigma(T_{|M_1})=\sigma_1$, which is not possible. Now consider the projection $P_1:X\rightarrow M_1$ along $M_1$ onto $M_2$. Since $J_{T}(x)\subset J_{T_{|M_1}}(x_1)\oplus J_{T_{|M_2}}(x_2)$, it follows that
\[P_1(J_{T}(x))\subset J_{T_{|M_1}}(x_1).\]

By the open mapping theorem it follows that $P_1(J_{T}(x))$ has non empty interior   and so $J_{T_{|M_1}}(x_1)$ has non-empty interior.
This is only possible if $x_1=0$, since $M_1$ is finite dimensional, see \cite{G.C.}.\\
Like above we conclude $J_{T_{|M_2}}(x_2)$ has non-empty interior. Therefore since $\sigma(T_{|M_2})=\sigma_2\subset\mathbb{C}\backslash\overline{\mathbb{D}}$, it follows by Lem. 1.8 that $x_2=0$. So alltogether $x=x_1+x_2=0$. This contradicts that $x\neq0$.\\
Now we proof the second statement.\\
If $\sigma(T)\subset\mathbb{C}\backslash\overline{\mathbb{D}}$ it is well known that $T$ is zero J-class, see \cite{G.C.}.
Now suppose $T$ is zero J-class and $\sigma(T)$ is not contained in $\mathbb{C}\backslash\overline{\mathbb{D}}$. Then like above we can decompose $X=M_1\oplus M_2$ and $T=T_1\oplus T_2$, where $M_1$ is finite dimensional. Since $T_1=T_{|M_1}$ is also zero J-class it follows, that $\sigma(T_1)\subset\mathbb{C}\backslash\overline{\mathbb{D}}$ (see \cite{G.C.}), which is a contradiction.
\end{proof}

Now S.Argyros, A.Arvanitakis and A.Tolias constructed a non-separable real Banach space, on which every operator $T$ has the form $T=\lambda I+S$, where $S$ is strictly singular and has separable range, see \cite{Argyros}.
\begin{theorem}\textbf{(Argyros, Arvanitakis, Tolias)}
There exists a real non-separable Banach space $X_{A}$ on which every operator $T$ is of the form $T=\lambda I+S$ $(\lambda\in\mathbb{R})$, where $S$ is strictly singular and has separable range.
\end{theorem} 

\begin{corollary} Consider $X:=(X_{A})_{\mathbb{C}}$. Then every operator $T\in L(X)$ is of the form $T=wI+S$ $(w\in\mathbb{C})$, where $S$ is strictly singular and has separable range.
\end{corollary}
\begin{proof} Every operator $T\in L(X)$ can be written as $T=T_1+iT_2$, where $T_1,T_2\in L(X_{A})$. By the previous theorem $T_1=\lambda I+S_1$ and $T_2=\mu I+S_2$ $(\lambda,\mu\in\mathbb{R})$, where $S_1,S_2\in L(X_A)$ are strictly singular and have separable range. Therefore we get
\begin{align*}
 T&=T_1+T_2\\
  &=\lambda I+S_1+i(\mu I+S_2)\\
  &= (\lambda+i\mu)I_{\mathbb{C}}+(S_1)_{\mathbb{C}}+i(S_2)_{\mathbb{C}}
 \end{align*}
Now by Prop. 1.7 $(S_i)_{\mathbb{C}}$ is strictly singular for $i\in\{1,2\}$ and by Thm. 1.6 $S:=(S_1)_{\mathbb{C}}+i(S_2)_{\mathbb{C}}$ is strictly singular and has separable range. With $w:=\lambda+i\mu$ we get $T=wI+S$.
\end{proof}

The next lemma can be found in \cite{G.C.}.
\begin{lemma} Let $X$ be a Banach space and $T\in L(X)$. If $J_T(x)$ has nonempty interior for some $x\neq 0$, then $T-\lambda I$ has dense range for each $\left|\lambda\right|\leq 1$.
\end{lemma}

\begin{theorem}
There exist a non-separable complex Banach space $X$ on which the J-set of every operator has empty interior for every non-zero vector. In particular there does not exist a J-class operator on $X$.
\end{theorem}
\begin{proof}
We condsider the space $X=(X_A)_{\mathbb{C}}$. Then every operator $T$ is of the form $T=\lambda I+S$ by Cor. 1.11, where $S$ is strictly singular and has separable range. If $\left|\lambda \right|>1$, then it follows from Thm. 1.9, that the interior of $J_{T}(x)$ is empty for each non-zero vector $x$. Now consider $\left|\lambda\right|\leq1$. Then by Lem. 1.12 the operator
$T-\lambda id=S$ has dense range. This is not possible since $S$ has separable range
and $X$ is non-separable.
\end{proof}
Our next aim is to show that on the space $Y:=X\oplus X$, where $X=(X_A)_{\mathbb{C}}$, the J-set of every $T\in L(Y)$ has also emtpty interior for each non-zero vector in $Y$. The next Lemma gives us some information about the form of the operators in $L(Y)$.
\begin{lemma} Consider $Y:=X\oplus X$, where $X=(X_A)_{\mathbb{C}}$. Then for every operator $T\in L(Y)$ there exists an isomorphism $J\in L(Y)$, such that $J^{-1}TJ$ has the following two possible matrix representations:
\[J^{-1}TJ=\begin{pmatrix} \lambda_1 I &  I\\ 0 & \lambda_2 I  \end{pmatrix}+\begin{pmatrix}S_1&S_2\\S_3&S_4\end{pmatrix}\]
or
\[J^{-1}TJ=\begin{pmatrix} \lambda_1 I &  0\\ 0 &  \lambda_2 I \end{pmatrix}+\begin{pmatrix}S_1&S_2\\S_3&S_4\end{pmatrix},\]
where $S_i\in L(X)$ is strictly singular for $i\in\{1,2,3,4\}$. 
\end{lemma}
\begin{proof}
Every operator $T\in L(Y)$ has the following matrix representation:
\[T=\begin{pmatrix} T_1 & T_2\\T_3 & T_4\end{pmatrix},\]
where $T_i\in L(X)$ for $i\in\{1,2,3,4\}$. Now by Cor 1.11 every $T_i=a_iI+\widetilde{S_i}$ with $\widetilde{S_i}$ strictly singular, so we get
\[T=\overbrace{\begin{pmatrix}a_1I&a_2I\\a_3I&a_4I\end{pmatrix}}^{A}+\overbrace{\begin{pmatrix}\widetilde{S_1}&\widetilde{S_2}\\ \widetilde{S_3}&\widetilde{S_4}\end{pmatrix}}^{\widetilde{S}}\]
Similiar to the Jordan decomposition of matrices, there exists now an isomorphism such that
\[J^{-1}AJ=\begin{pmatrix} \lambda_1 I &  I\\ 0 & \lambda_2 I  \end{pmatrix}\text{ or }J^{-1}AJ=\begin{pmatrix} \lambda_1 I &  0\\ 0 & \lambda_2 I  \end{pmatrix}\]
By Thm. 1.6 $S:=J^{-1}\widetilde{S}J$ is also strictly singular and therefore there exists some $S_i\in L(X)$, $i\in\{1,2,3,4\}$ strictly singular (see \cite{Abramovich}), such that
\[S=\begin{pmatrix}S_1&S_2\\S_3&S_4\end{pmatrix}.\]
It now follows the desired statement.
\end{proof}
\begin{theorem}
Consider $Y=X\oplus X$ with $X=(X_A)_{\mathbb{C}}$ and $T\in L(Y)$. Then $J_{T}((x,y))$ has empty interior for each non-zero vector $(x,y)\in Y$.
\end{theorem}
\begin{proof} 
We argue by contradiction. So suppose there exist some  $T\in L(Y)$, such that $\stackrel{\circ}{J_{T}((x,y))}\neq\emptyset$ for some $(x,y)\in Y\backslash\{0\}$. By the above lemma we will find an isomorphism $D\in L(Y)$ such that 
\[D^{-1}TD=\begin{pmatrix} \lambda_1 I &  I\\ 0 & \lambda_2 I  \end{pmatrix}+\begin{pmatrix}S_1&S_2\\S_3&S_4\end{pmatrix}\quad(*)\]
or
\[D^{-1}TD=\begin{pmatrix} \lambda_1 I &  0\\ 0 &  \lambda_2 I \end{pmatrix}+\begin{pmatrix}S_1&S_2\\S_3&S_4\end{pmatrix}\quad(**).\]
Then for $\widetilde{T}:=D^{-1}TD$ the J-set $J_{\widetilde{T}}(D^{-1}(x,y))$ has also non-empty interior. \\\\
$\textbf{Case 1}: \widetilde T=(*)\\ $
In this case $\lambda=\lambda_1=\lambda_2$. If $\left|\lambda\right|\neq1$ then we decompose $\sigma(\widetilde{T})$ in $\sigma_1=\{\mu\in\sigma(\widetilde{T}):\left|\mu\right|=1\}$ and $\sigma_2=\{\mu\in\sigma(\widetilde{T}):\left|\mu\right|\neq1\}$.
$\sigma_1$ is closed and by Thm. 1.4 $\sigma_2$ is also closed, since $\lambda$ and $0$ are the only possible accumulation points of $\sigma(\widetilde{T})$ and hence of $\sigma_2$. Furthermore the corresponding $\widetilde{T}$-invariant closed subspace $M_1$ for $\sigma_1$ is finite dimensional, otherwise 
\[\widetilde{T}|_{M_1}-\begin{pmatrix} \lambda I|_{M_1} &  0\\ 0 & \lambda I|_{M_1}  \end{pmatrix}=\begin{pmatrix}S_1|_{M_1}& I|_{M_1} +S_2|_{M_1}\\S_3|_{M_1}&S_4|_{M_1}\end{pmatrix}\]
is not invertible and hence $\lambda\in\sigma (\widetilde{T}|_{M_1})=\sigma_1$, which is not possible. The rest of the proof for $\left|\lambda\right|>1$ is analogous to Thm. 1.9\\
Now condsider $\left|\lambda\right|=1$. Then $\widetilde{T}-\begin{pmatrix}\lambda I&0\\0&\lambda I\end{pmatrix}$ has not dense range, which is a contradiction to lemma 1.12\\
\textbf{Case 2}:\ $\widetilde{T}=(**)$\\
If $\lambda_1=\lambda_2$ the argumentation is almost indentical to case 1.
So suppose $\lambda_1\neq\lambda_2$. Assume  $\left|\lambda_1\right|=1$ or $\left|\lambda_2\right|=1$. Without loss of generality $\left|\lambda_1\right|=1$. Then $\widetilde {T}-\begin{pmatrix}\lambda_1 I&0\\0& \lambda_1 I\end{pmatrix}$ has not dense range, which is a contradiction as in Case 1.\\
And for $\left|\lambda_1\right|\neq1$ and $\left|\lambda_2\right|\neq 1$ the argumentation is indentical like in Case 1. 
\end{proof}
\begin{remark}
It is also possible with some more technicalities to prove the same result in Thm 1.15 for $Y=\overbrace{X\oplus\ldots\oplus X}^{n-times}$, where $X=(X_{A})_{\mathbb{C}}$.
\end{remark}

We will now show that there is a large class of non-separable Banach spaces on which there exists always a J-class operator, namely the reflexive non-separable Banach spaces. The next theorem due to Costakis and Manoussos can be found in \cite{G.C.}.

\begin{theorem}
Let $X$ be a Banach space and $Y$ a separable Banach space. Consider $S\in L(X)$ with $\sigma(S)\subset\{\lambda:\left|\lambda\right|>1\}$. Let also $T\in L(Y)$ be hypercyclic. Then
\begin{enumerate}
\item $S\times T: X\times Y\rightarrow X\times Y$ is a J-class operator, but not hypercyclic.
\item $A_{S\times T}=\{0\}\times{Y}$.
\end{enumerate}
\end{theorem}

The next theorem by Lindenstrauss (\cite{J.L.}) gives us some information about the decomposition of reflexive non-separable Banach spaces. 
\begin{theorem}\textbf{(Lindenstrauss)}
Let $X$ be a non-separable reflexive Banach space and $Y\subset X$ a separable and closed subspace. Then there exists a separable closed subspace $W$ of $X$, which contains $Y$  and a linear bounded projection $P_W:X\rightarrow W$ with $\left\|P_W\right\|=1$.
\end{theorem}
\begin{theorem}
Let $X$ be a non-separable reflexive Banach space. Then for every infinite dimensional separable and closed subspace $Y$ and for every $\lambda\in (1,\infty)$ there exists a J-class operator $T$ with $Y\subset A_{T}$ and $\left\|T\right\|=\lambda$.
\end{theorem}
\begin{proof}
By Thm 1.19 there exists a separable infinite dimensional subspace $W$, which contains $Y$ and a linear bounded projection $P_{W}:X\rightarrow W$ with $\left\|P_{W}\right\|=1$.
There exists now a closed subspace $U$ of $X$ such that $X=U\oplus W$. For given $\epsilon>0$ we can find a hypercyclic operator $T_1:W\rightarrow W$, $T_1:= I_{W}+K$, with $K$ compact and $\left\|K\right\|<\epsilon$, see \cite{S.A.}. Then by Thm. 1.17 the operator $T_{\lambda}:=\lambda I_{U}\oplus T_1=\lambda I+(1-\lambda)P_{W}+K\circ P_{W}$ is J-class for $\lambda>1$. Furthermore $Y\subset W=A_{T}$.
Now define the function $g:(1,\infty)\rightarrow \mathbb{R}$ by $g(\delta):=\left\|T_{\delta}\right\|$. Then it is easy to see, that $g$ is continuous.
For given $\lambda$ we choose $\delta>1$ and $\epsilon>0$, such that $2\delta+\epsilon<1+\lambda$. Therefore we get
\[g(\delta)=\left\|T_{\delta}\right\|=\left\|\delta I+(1-\delta)P_{W}+K\circ P_{W}\right\| \leq\] \[\delta+\left|1-\delta\right|\left\|P_{W}\right\|+\left\|P_{W}\left\|\right\|K\right\|\leq 2\delta-1+\epsilon<\lambda.\]
On the other hand we can find a $\mu>1$ large enough, such that $g(\mu)>\lambda$. By the intermediate value theorem there exist now a $\xi\in\left[\delta,\mu\right]$ with $g(\xi)=\left\|T_{\xi}\right\|=\lambda$.
\end{proof}

\end{document}